\documentclass[12pt,leqno,draft]{article}

\newtheorem{theorem}{Theorem}
\newtheorem{lemma}[theorem]{Lemma}
\newtheorem{proposition}[theorem]{Proposition}
\newtheorem{definition}[theorem]{Definition}
\newtheorem{corollary}[theorem]{Corollary}

\newcommand{\begintheorem}{\addtocounter{equation}{1}\begin{theorem}}
\newcommand{\beginlemma}{\addtocounter{equation}{1}\begin{lemma}}
\newcommand{\beginproposition}{\addtocounter{equation}{1}\begin{proposition}}
\newcommand{\begindefinition}{\addtocounter{equation}{1}\begin{definition}}
\newcommand{\begincorollary}{\addtocounter{equation}{1}\begin{corollary}}

\begin{document}

\title{Some remarks about Cantor sets}

\author{Stephen William Semmes	\\
	Rice University		\\
	Houston, Texas}

\date{}

\maketitle

	The classical middle-thirds Cantor set can be described as
follows.  Start by taking $C_0 = [0, 1]$, the unit interval in the
real line.  Then put
\begin{equation}
	C_1 = [0, 1/3] \cup [2/3, 1],
\end{equation}
which is to say that one removes the open middle third from the unit
interval to get a union of two disjoint closed intervals of length
$1/3$.  By repeating the process one gets for each nonnegative integer
$j$ a subset $C_j$ of the unit interval which is a union of $2^j$
disjoint closed intervals of length $3^{-j}$ and which satisfies
\begin{equation}
	C_{j+1} \subseteq C_j
\end{equation}
for all $j$.  The Cantor set is then defined by
\begin{equation}
	C = \bigcap_{j=1}^\infty C_j.
\end{equation}

	For a more symbolic description it is convenient to consider
the set $\mathcal{B}$ of all binary sequences, which is to say all
sequences $x = \{x_l\}_{l=1}^\infty$ such that each term $x_l$ is
either equal to $0$ or $1$.  If $x$ is an element of $\mathcal{B}$,
consider
\begin{equation}
	\tau(x) = \sum_{l=1}^\infty \frac{2 x_l}{3^l}.
\end{equation}
Thus $\tau(x)$ is the real number in the unit interval given by a
ternary, or base $3$, expansion whose terms are $2 x_l$, and in
particular each term is equal to $0$ or $2$, while general ternary
expansions use $0$, $1$, and $2$.  It is well known that $\tau(x) \in
C$ for all $x \in \mathcal{B}$, and in fact that $\tau$ defines a
one-to-one correspondence between $\mathcal{B}$ and the Cantor set
$C$.

	Of course there is another well-known correspondence between
binary sequences and real numbers, which associates to $x \in
\mathcal{B}$ the real number
\begin{equation}
	\beta(x) = \sum_{l=1}^\infty \frac{x_l}{2^l}.
\end{equation}
Thus $\beta(x) \in [0, 1]$ and $\beta(x)$ has binary expansion given
by $x$.  Every real number in the unit interval arises in this manner,
and the binary expansion is unique except for a countable set of real
numbers for which there are two expansions.  The exceptions occur
for real numbers which can be expressed by binary expansions which are
eventually constant.

	More generally, suppose that for each positive integer $l$ we
have chosen a finite set $E_l$.  Define $\mathcal{E}$ to be the set of
sequences $x = \{x_l\}_{l=1}^\infty$ such that $x_l \in E_l$ for all
$l$.  For each $x \in \mathcal{E}$ and nonnegative integer $m$, define
$N_m(x)$ to be the subset of $\mathcal{E}$ consisting of sequences $y
= \{y_l\}_{l=1}^\infty$ such that $y_l = x_l$ when $l \le m$.
We can make $\mathcal{E}$ into a topological space by saying that
a subset $U$ of $\mathcal{E}$ is open if for each $x \in U$ there
is a nonnegative integer $m$ such that
\begin{equation}
	N_m(x) \subseteq U.
\end{equation}

	It is easy to see that $\mathcal{E}$ is a Hausdorff
topological space, and it is well known that $\mathcal{E}$ is also
compact.  It can be convenient to look at this in terms of sequential
compactness, which is to say that every sequence in $\mathcal{E}$
has a subsequence which converges.  It is easy to define metrics
on $\mathcal{E}$ for which the standard neighborhoods $N_m(x)$
are the balls.  Indeed, one can do this with ultrametrics, which
are metrics $d(x, y)$ which satisfy
\begin{equation}
	d(x, z) \le \max(d(x, y), d(y, z)).
\end{equation}

	In particular, one can define a topology on the set
$\mathcal{B}$ of binary sequences in this way.  The correspondence
$\tau$ between $\mathcal{B}$ and the Cantor set $C$ is then a
homeomorphism.  The mapping $\beta$ from $\mathcal{B}$ onto the unit
interval is a continuous mapping.  These are instances of a classical
and broad idea, which is that spaces can often be naturally
parameterized by spaces of the form $\mathcal{E}$.

	As a variant of this, suppose that for each integer $l$ we
have a finite set $E_l$, and that we specify a basepoint $b_l \in
E_l$.  Define $\mathcal{E}_*$ to be the set of doubly-infinite
sequences $x = \{x_l\}_{l = -\infty}^\infty$ such that $x_l \in E_l$
for all $l$ and there is an integer $L(x)$ with the property that $x_l
= b_l$ when $l \le L(x)$.  We can define standard neighborhoods
$N_m(x)$ in $\mathcal{E}_*$ for all $x \in \mathcal{E}_*$ and all
integers $m$ in exactly the same manner as before.  This leads to a
topology on $\mathcal{E}_*$ it is a locally compact Hausdorff space.

	These types of constructions show up in logic and computer
science as in \cite{1, 3}.  In another direction, suppose that each
$E_l$ is a finite semigroup, which is to say that it is a finite set
equipped with an associative binary operation.  One might also require
that each $E_l$ be equipped with an identity element, which can serve
as a basepoint in the second setting.  In this way $\mathcal{E}$ or
$\mathcal{E}_*$, as appropriate, can also be a semigroup, defining
a binary operation on sequences by applying the individual binary
operations to each term.  In fact one gets a topological semigroup,
which is to say that the binary operation is continuous.

	If the $E_l$'s are all groups, then one gets a group, with the
inverse operation being continuous as well.  There are other ways in
which topological groups have a similar topological structure, without
simply being a product.  A basic case is that of local fields, such as
the $p$-adic numbers \cite{2}, which are fields that are also locally
compact Hausdorff spaces in which the field operations are continuous,
and which are totally disconnected, in contrast with the real and
complex numbers.  Let us also mention profinite groups, which are
compact groups that are approximately finite in terms of a rich supply
of continuous homomorphisms to finite groups, which can also be
described in terms of closed subgroups of finite index.

	In symbolic dynamics one starts with a finite set $E$ with at
least two elements and considers the space $\mathcal{E}$ of all
sequences with terms in $E$, or the space $\widehat{\mathcal{E}}$ of
doubly-infinite sequences with terms in $E$.  For the latter we can
consider neighborhoods of a point $x \in \widehat{\mathcal{E}}$
consisting of all $y \in \widehat{\mathcal{E}}$ such that
\begin{equation}
	y_l = x_l \hbox{ when } |l| \le m
\end{equation}
for some nonnegative integer $m$.  This leads to a topology on
$\widehat{\mathcal{E}}$ in which it becomes a compact Hausdorff
topological space.  As a topological space, $\mathcal{E}$ and
$\widehat{\mathcal{E}}$ use essentially the same construction, except
in terms of formatting.

	In both situations there is a natural shift mapping which
takes a sequence $\{x_l\}$ to the shifted sequence $\{x_{l+1}\}$.  In
the case of $\mathcal{E}$ this is a continuous mapping of
$\mathcal{E}$ onto itself which is not one-to-one because the first
term in the sequence is dropped.  In the case of
$\widehat{\mathcal{E}}$ this defines a homeomorphism from
$\widehat{\mathcal{E}}$ onto itself.

	It is also frequently natural to look at measures and
integration on sets like this.  Let us start with the case of finitely
many finite sets $E_1, \ldots, E_n$, and with $\mathcal{E}$ defined to
be the set of $n$-tuples $(x_1, \ldots, x_n)$ such that $x_i \in E_i$
for $i = 1, \ldots, n$.  By a measure on $\mathcal{E}$ we simply mean
an assignment of a nonnegative real number to each element of
$\mathcal{E}$.  In many situations one might wish to require that the
measure of each point in $\mathcal{E}$ is positive.  One might also
wish to ask that the total measure of $\mathcal{E}$ be equal to $1$.

	A very important way in which such measures arise is from
products of measures on the individual $E_i$'s.  It may be that the
$E_i$'s are all the same, in which case one might start with a
measure on a single set $E$ and get an induced product measure
on the Cartesian product of $n$ copies of $E$.  At any rate,
it may also be that on each $E_i$ one uses a uniform distribution,
in which the measure of each element of $E_i$ is the same.

	Now suppose that we have chosen a finite set $E_l$ for each
positive integer $l$, and we let $\mathcal{E}$ be the space of
sequences with terms in the $E_l$'s as described previously.  One way
to look at measures is in terms of measures of the standard
neighborhoods $N_m(x)$, with the condition that for a finite
collection of such neighborhoods whose union is a standard
neighborhood the sum of the measures of the initial neighborhoods is
equal to the measure of their union.  This basically amounts
to having measures on the finite products
\begin{equation}
	E_1 \times \cdots \times E_n
\end{equation}
with suitable compatibility conditions as $n$ increases.

	Another way to look at this is to have a linear mapping from
continuous complex-valued functions on $\mathcal{E}$ into the complex
numbers which is nonnegative in the sense that the ``integral'' of a
nonnegative real-valued continuous function on $\mathcal{E}$ is a
nonnegative real number.  In practice this might be defined like a
classical Riemann integral, namely, as a limit of finite sums, where
the sums are based on measures of standard neighborhoods.  Conversely,
if one starts with an integral of continuous functions like this, then
one gets a measure for standard neighborhoods by integrating the
functions which are equal to $1$ on a standard neighborhood and which
are equal to $0$ on the rest of $\mathcal{E}$.  At any rate, this is a
rather friendly aspect of generalized Cantor sets like this, that
there are plenty of continuous functions which are locally constant,
and that arbitrary continuous functions can be approximated by
locally-constant functions.

\end{document}